# FROM FINITE SAMPLE TO ASYMPTOTICS: A GEOMETRIC BRIDGE FOR SELECTION CRITERIA IN SPLINE REGRESSION[1]

By S. C. Kou

*Harvard University*

This paper studies, under the setting of spline regression, the connection between finite-sample properties of selection criteria and their asymptotic counterparts, focusing on bridging the gap between the two. We introduce a bias-variance decomposition of the prediction error, using which it is shown that in the asymptotics the bias term dominates the variability term, providing an explanation of the gap. A geometric exposition is provided for intuitive understanding. The theoretical and geometric results are illustrated through a numerical example.

**1. Introduction.** A central problem in statistics is regression: One observes $\{(x_i, y_i), i = 1, 2, \ldots, n\}$ and wants to estimate the regression function of $y$ on $x$. Through the efforts of many authors, the past two decades have witnessed the establishment of nonparametric regression as a powerful tool for data analysis; references include, for example, Härdle (1990), Hastie and Tibshirani (1990), Wahba (1990), Silverman (1985), Rosenblatt (1991), Green and Silverman (1994), Eubank (1988), Simonoff (1996), Fan and Gijbels (1996), Bowman and Azzalini (1997) and Fan (2000).

The practical application of nonparametric regression typically requires the specification of a smoothing parameter which crucially determines how locally the smoothing is done. This article, under the setting of smoothing splines, concerns the data-driven choice of smoothing parameter (as opposed to a subjective selection); in particular, this article focuses on the connection between *finite-sample* properties of selection criteria and their *asymptotic* counterparts.

Received March 2003; revised March 2004.
[1]Supported in part by NSF Grant DMS-02-04674 and Harvard University Clark-Cooke Fund.

The *large-sample* (asymptotic) perspective has been impressively addressed in the literature. Some references, among others, include Wahba (1985), Li (1986, 1987), Stein (1990), Hall and Johnstone (1992), Jones, Marron and Sheather (1996), Hurvich, Simonoff and Tsai (1998) and Speckman and Sun (2001).

Complementary to the large-sample (asymptotic) developments, Efron (2001) and Kou and Efron (2002), using a geometric interpretation of selection criteria, study the *finite-sample* properties. For example, they explain (a) why the popular $C_p$ criterion has the tendency to be highly variable [even for data sets generated from the same underlying curve, the $C_p$-estimated curve varies a lot from oversmoothed ones to very wiggly ones; see Kohn, Ansley and Tharm (1991) and Hurvich, Simonoff and Tsai (1998) for examples], and (b) why another selection criterion, generalized maximum likelihood [Wecker and Ansley (1983), Wahba (1985) and Stein (1990)], appears to be stable and yet sometimes tends to undersmooth the curve. Roughly speaking, it was shown that the root of the variable behavior of $C_p$ is its geometric instability, while the stable but undersmoothing behavior of generalized maximum likelihood (GML) stems from its potentially large bias. In addition, they also introduce a new selection criterion, the extended exponential (EE) criterion, which combines the strength of $C_p$ and GML while mitigating their weaknesses.

With the asymptotic and finite-sample properties delineated, it seems that we have a "complete" picture of selection criteria. However, a careful inspection of the finite-sample and asymptotic results, especially the ones comparing $C_p$ and GML, reveals an interesting gap. On the finite-sample side, $C_p$'s geometric instability undermines its competitiveness [Kohn, Ansley and Tharm (1991) and Hurvich, Simonoff and Tsai (1998)], which opens the door for the more stable GML, while on the large-sample (asymptotic) side different authors [e.g., Wahba (1985) and Li (1986, 1987)] have suggested that from the frequentist standpoint the $C_p$-type criterion asymptotically performs more efficiently than GML. This "gap" between finite-sample and asymptotic results naturally makes one puzzle: (a) Why doesn't the finite-sample advantage of GML, notably its stability, benefit it as far as large-sample (asymptotics) is concerned? (b) Why does the geometric instability of $C_p$ seen in finite-sample disappear in the asymptotic considerations?

This article attempts to address these puzzles. First, by decomposing the estimation error into a bias part and a variability part, we show that as sample size grows large the bias term dominates the variability term, thus making the large-sample case virtually a bias problem. Consequently in the large-sample comparisons, one is essentially comparing the bias of different selection criteria and unintentionally overlooking the variability—a situation particularly favoring the $C_p$-type criterion as it is (asymptotically) unbiased. Second, by studying the evolution of the geometry of selection criteria, we show that the geometric instability of selection criteria gradually decreases, though rather slowly, which again benefits the $C_p$-type criterion, because it says as far as asymptotics is concerned, the instability of $C_p$ evident in finite-sample studies will not show up. The recent interesting work of Speckman and Sun (2001) appears to confirm our results regarding asymptotics (see Section 2); they showed that GML and $C_p$ agree on the relative convergence rate of the selected smoothing parameter.

The connection between finite-sample and asymptotic results is illustrated by a numerical example (Section 4). The numerical example also indicates that for sample sizes one usually encounters in practice the EE criterion appears to behave more stably than both GML and $C_p$.

The article is organized as follows. Section 2 introduces a bias-variance decomposition of the total prediction error, and investigates its finite- and large-sample consequences. Section 3 provides a geometric explanation to bridge the finite-sample and asymptotic results regarding selection criteria. Section 4 illustrates the connection through a simulation experiment. The article concludes in Section 5 with further remarks. The detailed theoretical proofs are deferred to the Appendix.

## 2. A bias-variance decomposition for prediction error.

2.1. *Selection criteria in spline regression.* The goal of regression is to estimate $f(x) = E(y|x)$ from $n$ observed data points $\{(x_i, y_i), i = 1, 2, \ldots, n\}$. A linear smoother estimates $\mathbf{f} = (f(x_1), f(x_2), \ldots f(x_n))'$ by $\hat{\mathbf{f}}_\lambda = \mathbf{A}_\lambda \mathbf{y}$, where the entries of the $n \times n$ smoothing matrix $\mathbf{A}_\lambda$ depend on $\mathbf{x} = (x_1, x_2, \ldots, x_n)$ and also on a nonnegative *smoothing parameter* $\lambda$. One class of linear smoothers that will be of particular interest in this article is (cubic) smoothing splines, under which

$$(2.1) \qquad \mathbf{A}_\lambda = \mathbf{U}\mathbf{a}_\lambda \mathbf{U}',$$

where $\mathbf{U}$ is an $n \times n$ orthogonal matrix *not* depending on $\lambda$, and $\mathbf{a}_\lambda = \mathrm{diag}(a_{\lambda i})$, a diagonal matrix with the $i$th diagonal element $a_{\lambda i} = 1/(1 + \lambda k_i)$, $i = 1, 2, \ldots, n$. The constants $\mathbf{k} = (k_1, k_2, \ldots, k_n)$, solely determined by $\mathbf{x}$, are nonnegative and nondecreasing. The trace of the smoothing matrix $\mathrm{tr}(\mathbf{A}_\lambda)$ is referred to as the "degrees of freedom," $df_\lambda = \mathrm{tr}(\mathbf{A}_\lambda)$, which agrees with the standard definition if $\mathbf{A}_\lambda$ represents polynomial regression.

To use splines in practice, one typically has to infer the value of the smoothing parameter $\lambda$ from the data. The $C_p$ criterion chooses $\lambda$ to minimize an unbiased estimate of the total squared error. Suppose the $y_i$'s are uncorrelated, with mean $f_i$ and constant variance $\sigma^2$. The $C_p$ estimate of $\lambda$ is $\hat{\lambda}^{C_p} = \arg\min_\lambda \{C_\lambda(\mathbf{y})\}$, where the $C_p$ statistic $C_\lambda(\mathbf{y}) = \|\mathbf{y} - \hat{\mathbf{f}}_\lambda\|^2 + 2\sigma^2 \mathrm{tr}(\mathbf{A}_\lambda) - n\sigma^2$ is an unbiased estimate of $E\|\hat{\mathbf{f}}_\lambda - \mathbf{f}\|^2$.

The generalized maximum likelihood (GML) criterion [Wecker and Ansley (1983)] is another selection criterion motivated from empirical Bayes considerations. If one starts from $\mathbf{y} \sim N(\mathbf{f}, \sigma^2 \mathbf{I})$, and puts a Gaussian prior on the underlying curve: $\mathbf{f} \sim N(\mathbf{0}, \sigma^2 \mathbf{A}_\lambda (\mathbf{I} - \mathbf{A}_\lambda)^{-1})$, then by Bayes theorem,

$$(2.2) \qquad \mathbf{y} \sim N(\mathbf{0}, \sigma^2 (\mathbf{I} - \mathbf{A}_\lambda)^{-1}), \qquad \mathbf{f}|\mathbf{y} \sim N(\mathbf{A}_\lambda \mathbf{y}, \sigma^2 \mathbf{A}_\lambda).$$

The second relationship shows that $\hat{\mathbf{f}}_\lambda = \mathbf{A}_\lambda \mathbf{y}$ is the Bayes estimate of $\mathbf{f}$. The first relationship motivates the GML: It chooses $\hat{\lambda}^{\mathrm{GML}}$ as the MLE of $\lambda$ from $\mathbf{y} \sim N(\mathbf{0}, \sigma^2 (\mathbf{I} - \mathbf{A}_\lambda)^{-1})$.

The setting of smoothing splines (2.1) allows a rotation of coordinates,

$$(2.3) \qquad \mathbf{z} = \mathbf{U}' \mathbf{y}/\sigma, \qquad \mathbf{g} = \mathbf{U}' \mathbf{f}/\sigma, \qquad \hat{\mathbf{g}}_\lambda = \mathbf{U}' \hat{\mathbf{f}}_\lambda /\sigma,$$

which leads to a diagonal form: $\mathbf{z} \sim N(\mathbf{g}, \mathbf{I})$, $\hat{\mathbf{g}}_\lambda = \mathbf{a}_\lambda \mathbf{z}$. Let $b_{\lambda i} = 1 - a_{\lambda i}$ and $\mathbf{b}_\lambda = (b_{\lambda 1}, b_{\lambda 2}, \ldots, b_{\lambda n})$. In the new coordinate system, the $C_p$ statistic can be expressed as a function of $\mathbf{z}^2$, $C_\lambda(\mathbf{z}^2) = \sigma^2 \sum_{i=1}^n (b_{\lambda i}^2 z_i^2 - 2 b_{\lambda i}) + n\sigma^2$, and correspondingly

$$\hat{\lambda}^{C_p} = \arg\min_\lambda \sum_{i>2} (b_{\lambda i}^2 z_i^2 - 2 b_{\lambda i}).$$

Under the coordinate system of $\mathbf{z}$ and $\mathbf{g}$, since $\mathbf{z} \sim N(\mathbf{0}, \mathrm{diag}(\mathbf{b}_\lambda^{-1}))$, $\mathbf{g}|\mathbf{z} \sim N(\mathbf{a}_\lambda \mathbf{z}, \mathbf{a}_\lambda)$,

$$\hat{\lambda}^{\mathrm{GML}} = \text{MLE of } \mathbf{z} \sim N(\mathbf{0}, \mathrm{diag}(\mathbf{b}_\lambda^{-1})) = \arg\min_\lambda \sum_{i>2} (b_{\lambda i} z_i^2 - \log b_{\lambda i}).$$

Because $\mathbf{z}$ and $\mathbf{g}$ offer simpler expressions, we will work on them instead of $\mathbf{y}$ and $\mathbf{f}$ whenever possible. The extended exponential (EE) selection criterion, studied in Kou and Efron (2002), provides a third way to choose the smoothing parameter. It is motivated by the idea of combining the strengths of $C_p$ and GML while mitigating their weaknesses, since in practice the $C_p$-selected smoothing parameter tends to be highly variable, whereas the GML criterion has a serious problem with bias (see Section 4 for an illustration). Expressed in terms of $\mathbf{z}$, the EE criterion selects the smoothing parameter $\lambda$ according to

$$\hat{\lambda}^{\mathrm{EE}} = \arg\min_\lambda \sum_{i>2} [C b_{\lambda i} z_i^{4/3} - 3 b_{\lambda i}^{1/3}],$$

where the constant $C = \frac{\sqrt{\pi}}{2^{2/3} \Gamma(7/6)} = 1.203$. Kou and Efron (2002) explained its construction from a geometric point of view and illustrated through a finite-sample nonasymptotic analysis that the EE criterion combines the strengths of $C_p$ and GML to a large extent.

An interesting fact about the three criteria ($C_p$, GML and EE) is that they share a unified structure. Let $p \geq 1$, $q \geq 1$ be two fixed constants. Define the function

$$(2.4) \quad l_\lambda^{(p,q)}(\mathbf{u}) = \begin{cases} \sum_i \left[ (c_q b_{\lambda i}^{1/q})^p u_i - \dfrac{p}{p-1}((c_q b_{\lambda i}^{1/q})^{p-1} - 1) \right], & \text{if } p > 1, \\ \sum_i (c_q b_{\lambda i}^{1/q} u_i - \log b_{\lambda i}^{1/q}), & \text{if } p = 1, \end{cases}$$

where $c_q = \dfrac{\sqrt{\pi}}{2^{1/q}\Gamma(1/2+1/q)}$, and a corresponding selection criterion

$$(2.5) \quad \hat{\lambda}^{(p,q)} = \arg\min_\lambda \{l_\lambda^{(p,q)}(\mathbf{z}^{2/q})\}.$$

Then it is easy to verify that (i) $l_\lambda^{(p,q)} \to l_\lambda^{(1,q)}$ as $p \to 1$; (ii) taking $p = 1$, $q = 1$ gives the GML criterion; $p = 2$, $q = 1$ gives the $C_p$ criterion; $p = q = \frac{3}{2}$ gives the EE criterion. The class (2.5), therefore, unites the three criteria in a continuous fashion. This much facilitates our theoretical development as it allows us to work on the general selection criterion $\hat{\lambda}^{(p,q)}$ and take $(p,q)$ to specific values to obtain corresponding results for EE, $C_p$ and GML.

2.2. *The unbiasedness of $C_p$.* To introduce the idea of bias-variance decomposition, we first note that for each selection criterion $\hat{\lambda}^{(p,q)}$ there are an associated central smoothing parameter $\lambda_c^{(p,q)}$ and central degrees of freedom $df_c^{(p,q)}$ obtained by applying the expectation operator on the selection criterion (2.5):

$$(2.6) \quad \lambda_c^{(p,q)} = \arg\min_\lambda E\{l_\lambda^{(p,q)}(\mathbf{z}^{2/q})\},$$

$$(2.7) \quad df_c^{(p,q)} = \operatorname{tr}(\mathbf{A}_{\lambda_c^{(p,q)}}).$$

Since (2.6) is the estimating-equation version of (2.5), from the general theory of estimating equations it can be seen that $\hat{\lambda}^{(p,q)}$ and $\widehat{df}^{(p,q)}$ are centered around $\lambda_c^{(p,q)}$ and $df_c^{(p,q)}$ in the sense that $\lambda_c^{(p,q)}$ and $df_c^{(p,q)}$ are the asymptotic means of $\hat{\lambda}^{(p,q)}$ and $\widehat{df}^{(p,q)}$. Thus $\lambda_c^{(p,q)}$ and $df_c^{(p,q)}$ index the central tendency of the selection criterion-$(p,q)$.

Next we introduce the *ideal smoothing parameter* $\lambda_0$ and the *ideal degrees of freedom* $df_0 = \operatorname{tr}(\mathbf{A}_{\lambda_0})$, which are intrinsically determined by the underlying curve and do not depend on the specific selection criterion one uses:

$$(2.8) \quad \lambda_0 = \arg\min_\lambda E_{\mathbf{f}} \|\hat{\mathbf{f}}_\lambda - \mathbf{f}\|^2 = \arg\min_\lambda E\|\hat{\mathbf{g}}_\lambda - \mathbf{g}\|^2.$$

The risk $E\|\hat{\mathbf{g}}_{\lambda_0} - \mathbf{g}\|^2$ associated with $\lambda_0$ represents the minimum risk one has to bear to estimate the underlying curve. Therefore, to compare the

performance of different selection criteria one can focus on the *extra* risk: $E\|\hat{\mathbf{g}}_{\hat{\lambda}} - \mathbf{g}\|^2 - E\|\hat{\mathbf{g}}_{\lambda_0} - \mathbf{g}\|^2$. See Wahba (1985), Härdle, Hall and Marron (1988), Hall and Johnstone (1992), Gu (1998) and Efron (2001) for more discussion.

Having introduced the necessary concepts, we state our first result, the unbiasedness of $C_p$.

THEOREM 2.1. *The central smoothing parameter $\lambda_c^{(2,1)}$ and degrees of freedom $df_c^{(2,1)}$ of $C_p$ correspond exactly to the ideal smoothing parameter and degrees of freedom*

$$\lambda_c^{(2,1)} = \lambda_0, \qquad df_c^{(2,1)} = df_0.$$

PROOF. First, from the definition (2.8) a straightforward expansion gives

$$\lambda_0 = \arg\min_\lambda \sum_i (b_{\lambda i}^2 (g_i^2 + 1) - 2b_{\lambda i}). \tag{2.9}$$

Next, for $C_p$ according to (2.6) its central smoothing parameter

$$\begin{aligned}
\lambda_c^{(2,1)} &= \arg\min_\lambda E\{l_\lambda^{(2,1)}(\mathbf{z}^2)\} = \arg\min_\lambda E\left\{\sum_i [b_{\lambda i}^2 z_i^2 - 2b_{\lambda i}]\right\} \\
&= \arg\min_\lambda \sum_i (b_{\lambda i}^2 (g_i^2 + 1) - 2b_{\lambda i}).
\end{aligned} \tag{2.10}$$

The proof is complete because (2.9) and (2.10) give identical expressions for $\lambda_0$ and $\lambda_c^{(2,1)}$. □

Since no other element from the selection criteria class (2.5) possesses this property of unbiasedness, the result of Theorem 2.1 gives $C_p$ an advantage over the others. As we shall see shortly, this advantage is the main factor that makes the asymptotic consideration favorable for $C_p$.

2.3. *The bias-variance decomposition.* The results developed so far work for all sample sizes. Next we turn our attention to the large-sample case. There is a large amount of literature addressing the large-sample properties of selection criteria. The well-cited asymptotic results [Wahba (1985) and Li (1986, 1987), among others] suggest that as far as large-sample is concerned the $C_p$-type criterion outperforms GML. This interestingly seems at odds with the well-known finite-sample results. For example, Kohn, Ansley and Tharm (1991) and Hurvich, Simonoff and Tsai (1998), among others, illustrate that finite-sample-wise the $C_p$ criterion has a strong tendency for high variability in the sense that even for data sets generated from the same underlying curve the $C_p$-estimated curves vary a great deal from oversmoothed ones to very wiggly ones, which contrasts with the stably performing GML. To understand why there is this gap between finite- and large-sample results, we will provide a bias-variance decomposition of the prediction error, based on which it will be seen that the major reason is that the large-sample consideration virtually only looks at the bias, as bias asymptotically dominates variability.

The central smoothing parameter and central degrees of freedom defined previously pave the way for the bias-variance decomposition. Consider the prediction error for estimating the curve $E\|\hat{\mathbf{f}}_{\hat{\lambda}(p,q)} - \mathbf{f}\|^2$, which is equal to $\sigma^2 E\|\hat{\mathbf{g}}_{\hat{\lambda}(p,q)} - \mathbf{g}\|^2$ according to (2.3). We can write

$$E\|\hat{\mathbf{g}}_{\hat{\lambda}(p,q)} - \mathbf{g}\|^2$$
$$= E\|(\hat{\mathbf{g}}_{\hat{\lambda}(p,q)} - \hat{\mathbf{g}}_{\lambda_c^{(p,q)}}) + (\hat{\mathbf{g}}_{\lambda_c^{(p,q)}} - \mathbf{g})\|^2$$
$$= E\|\hat{\mathbf{g}}_{\lambda_c^{(p,q)}} - \mathbf{g}\|^2 + 2E(\hat{\mathbf{g}}_{\lambda_c^{(p,q)}} - \mathbf{g})'(\hat{\mathbf{g}}_{\hat{\lambda}(p,q)} - \hat{\mathbf{g}}_{\lambda_c^{(p,q)}}) + E\|\hat{\mathbf{g}}_{\hat{\lambda}(p,q)} - \hat{\mathbf{g}}_{\lambda_c^{(p,q)}}\|^2.$$

Consequently, the extra risk beyond the unavoidable risk $E\|\hat{\mathbf{g}}_{\lambda_0} - \mathbf{g}\|^2$ can be written as

$$E\|\hat{\mathbf{g}}_{\hat{\lambda}(p,q)} - \mathbf{g}\|^2 - E\|\hat{\mathbf{g}}_{\lambda_0} - \mathbf{g}\|^2$$
(2.11)
$$= (E\|\hat{\mathbf{g}}_{\lambda_c^{(p,q)}} - \mathbf{g}\|^2 - E\|\hat{\mathbf{g}}_{\lambda_0} - \mathbf{g}\|^2)$$
$$+ 2E(\hat{\mathbf{g}}_{\lambda_c^{(p,q)}} - \mathbf{g})'(\hat{\mathbf{g}}_{\hat{\lambda}(p,q)} - \hat{\mathbf{g}}_{\lambda_c^{(p,q)}}) + E\|\hat{\mathbf{g}}_{\hat{\lambda}(p,q)} - \hat{\mathbf{g}}_{\lambda_c^{(p,q)}}\|^2.$$

This expression provides a bias-variance decomposition for the prediction error. The first term $E\|\hat{\mathbf{g}}_{\lambda_c^{(p,q)}} - \mathbf{g}\|^2 - E\|\hat{\mathbf{g}}_{\lambda_0} - \mathbf{g}\|^2$ can be viewed as the bias term—it captures the error of estimating the curve $\mathbf{g}$ beyond the unavoidable risk by using the central smoothing parameter $\lambda_c^{(p,q)}$, which measures the discrepancy between the central risk associated with $\hat{\lambda}^{(p,q)}$ and the ideal minimum risk; the third term $E\|\hat{\mathbf{g}}_{\hat{\lambda}(p,q)} - \hat{\mathbf{g}}_{\lambda_c^{(p,q)}}\|^2$ can be viewed as the variability term—it measures the variability of $\hat{\mathbf{g}}_{\hat{\lambda}(p,q)}$ from its "center" $\hat{\mathbf{g}}_{\lambda_c^{(p,q)}}$; the second term, the covariance, arises here due to the nature of adaptation (the smoothing parameter itself is also inferred from the data, in addition to estimating the curve).

Clearly, for any practical finite-sample problem, each term in (2.11) contributes to the squared prediction error. However, we shall show that as the sample size $n$ grows large the bias term gradually dominates the other two. To focus on the basic idea, without loss of generality, we assume the design points $(x_1, x_2, \ldots, x_n)$ are $n$ equally spaced points along the interval $[0,1]$. Section 5 will discuss the setting of general design points.

In what follows, to avoid cumbersome notation, we will write $\hat{\lambda}$ for $\hat{\lambda}^{(p,q)}$, $\widehat{df}$ for $\widehat{df}^{(p,q)}$, $\lambda_c$ for $\lambda_c^{(p,q)}$, $df_c$ for $df_c^{(p,q)}$, and so on. The full notation $\hat{\lambda}^{(p,q)}$, $\lambda_c^{(p,q)}$, $df_c^{(p,q)}$ will be used whenever potential confusion might arise. Consider the bias term $E\|\hat{\mathbf{g}}_{\lambda_c} - \mathbf{g}\|^2 - E\|\hat{\mathbf{g}}_{\lambda_0} - \mathbf{g}\|^2$ first:

$$E\|\hat{\mathbf{g}}_{\lambda_c} - \mathbf{g}\|^2 = E\|\mathbf{a}_{\lambda_c}\mathbf{z} - \mathbf{g}\|^2 = \sum_{i=1}^{n}(b_{\lambda_c i}^2 g_i^2 + a_{\lambda_c i}^2)$$

(2.12)

$$= \lambda_c \sum_{i=1}^{n}[a_{\lambda_c i}b_{\lambda_c i}(k_i g_i^2)] + \sum_{i=1}^{n} a_{\lambda_c i}^2,$$

where the last equality uses the fact $b_{\lambda i} = \frac{\lambda k_i}{1+\lambda k_i} = \lambda k_i a_{\lambda i}$. To obtain the asymptotic orders, we need to know how $\lambda_c$, the central smoothing parameter, evolves as the sample size gets large. According to definition (2.6), $\lambda_c$ satisfies the normal equation $\frac{\partial}{\partial \lambda} l_\lambda^{(p,q)}(E\{\mathbf{z}^{2/q}\})|_{\lambda=\lambda_c} = 0$, which (through some algebra) can be written as

(2.13) $$\sum_i a_{\lambda_c i} b_{\lambda_c i}^{p/q}(c_q E\{z_i^{2/q}\} - 1) = \sum_i a_{\lambda_c i} b_{\lambda_c i}^{(p-1)/q} - \sum_i a_{\lambda_c i} b_{\lambda_c i}^{p/q}.$$

The following lemma gives the order of the left-hand side of (2.13).

LEMMA 2.2. *Under mild regularity conditions, for $p \geq q$, $\sum_i a_{\lambda_c i} b_{\lambda_c i}^{p/q} \times (c_q E\{z_i^{2/q}\} - 1) = O(\lambda_c)$.*

The regularity conditions and the proof of Lemma 2.2 are given in the Appendix. The proof uses one handy result of Demmler and Reinsch (1975), where by studying the oscillation of the smoothing-spline eigenvectors, it is effectively shown that for any curve $f(x)$ satisfying $0 < \int_0^1 f''(t)^2 dt < \infty$,

(2.14) $$0 < \sum_{i=3}^{n} k_i g_i^2 \leq \frac{1}{\sigma^2} \int_0^1 f''(t)^2 dt < \infty \quad \text{for all } n \geq 3.$$

See also Speckman (1983, 1985) and Wahba (1985). For the right-hand side of (2.13), the following theorem, taken from Kou (2003), is useful.

THEOREM 2.3. *Suppose $\frac{n}{\lambda} \to \infty$ and $n^3 \lambda \to \infty$. Then for $r > \frac{1}{4}$, $s > -\frac{1}{4}$,*

$$\sum_{i=3}^{n} a_{\lambda i}^r b_{\lambda i}^s = \frac{1}{4\pi} B\left(r - \frac{1}{4}, s + \frac{1}{4}\right) \left(\frac{n}{\lambda}\right)^{1/4} + o\left(\left(\frac{n}{\lambda}\right)^{1/4}\right),$$

*where the beta function $B(x,y) = \Gamma(x)\Gamma(y)/\Gamma(x+y)$.*

Applying this result, the right-hand side of (2.13) is

$$\sum_i a_{\lambda_c i} b_{\lambda_c i}^{(p-1)/q} - \sum_i a_{\lambda_c i} b_{\lambda_c i}^{p/q} = O\left(\left(\frac{n}{\lambda_c}\right)^{1/4}\right).$$

Matching it with the result of Lemma 2.2 gives

(2.15) $$\lambda_c^{(p,q)} = O(n^{1/5}) \quad \text{for all } p \geq q,$$

which furthermore implies (taking $r = 1$, $s = 0$ in Theorem 2.3)

(2.16) $$df_c^{(p,q)} = O\left(\left(\frac{n}{\lambda_c^{(p,q)}}\right)^{1/4}\right) = O(n^{1/5}) \quad \text{for } p \geq q.$$

Note that (2.15) and (2.16) cover GML, $C_p$ and EE, since all three satisfy $p \geq q$. With the help of Theorem 2.3 and (2.15), we can calculate the asymptotic order of the bias term $E\|\hat{\mathbf{g}}_{\lambda_c} - \mathbf{g}\|^2 - E\|\hat{\mathbf{g}}_{\lambda_0} - \mathbf{g}\|^2$. By inequality (2.14), the first term of (2.12)

$$\lambda_c \sum_{i=1}^n [a_{\lambda_c i} b_{\lambda_c i}(k_i g_i^2)] \leq \lambda_c \sum_{i=1}^n (k_i g_i^2) = O(\lambda_c) = O(n^{1/5});$$

and (from Theorem 2.3) the second term of (2.12) $\sum_{i=1}^n a_{\lambda_c i}^2 = O((\frac{n}{\lambda_c})^{1/4}) = O(n^{1/5})$. Adding them together yields

(2.17) $$E\|\hat{\mathbf{g}}_{\lambda_c^{(p,q)}} - \mathbf{g}\|^2 = O(n^{1/5}).$$

Identical treatment of the ideal smoothing parameter $\lambda_0$ gives

(2.18) $$E\|\hat{\mathbf{g}}_{\lambda_0} - \mathbf{g}\|^2 = O(n^{1/5}).$$

Combining the results of (2.17) and (2.18), we observe that for a "general" criterion $\hat{\lambda}^{(p,q)}$ the bias term $E\|\hat{\mathbf{g}}_{\lambda_c^{(p,q)}} - \mathbf{g}\|^2 - E\|\hat{\mathbf{g}}_{\lambda_0} - \mathbf{g}\|^2 = O(n^{1/5})$. We put a quotation mark on "general" because there is one exception: $C_p$. In Theorem 2.1 we have shown that $\lambda_c^{(2,1)} = \lambda_0$, which implies that $E\|\hat{\mathbf{g}}_{\lambda_c^{(2,1)}} - \mathbf{g}\|^2 - E\|\hat{\mathbf{g}}_{\lambda_0} - \mathbf{g}\|^2 = 0$. The following theorem summarizes the discovery and extends the result to the variability and covariance terms in the decomposition.

THEOREM 2.4. *Under mild regularity conditions provided in the Appendix, for all $p \geq q$:*

(i) *the bias term*

$$E\|\hat{\mathbf{g}}_{\lambda_c^{(p,q)}} - \mathbf{g}\|^2 - E\|\hat{\mathbf{g}}_{\lambda_0} - \mathbf{g}\|^2 = \begin{cases} O(n^{1/5}), & \text{if } (p,q) \neq (2,1), \\ 0, & \text{if } (p,q) = (2,1), \end{cases}$$

(ii) *the covariance term* $E(\hat{\mathbf{g}}_{\lambda_c^{(p,q)}} - \mathbf{g})'(\hat{\mathbf{g}}_{\hat{\lambda}^{(p,q)}} - \hat{\mathbf{g}}_{\lambda_c^{(p,q)}}) = O(1)$,

(iii) *the variability term* $E\|\hat{\mathbf{g}}_{\hat{\lambda}^{(p,q)}} - \hat{\mathbf{g}}_{\lambda_c^{(p,q)}}\|^2 = O(1)$.

*Therefore, the extra risk*

$$E\|\hat{\mathbf{g}}_{\hat{\lambda}^{(p,q)}} - \mathbf{g}\|^2 - E\|\hat{\mathbf{g}}_{\lambda_0} - \mathbf{g}\|^2 = \begin{cases} O(n^{1/5}), & \text{if } (p,q) \neq (2,1), \\ O(1), & \text{if } (p,q) = (2,1). \end{cases}$$

The regularity conditions and the proof of Theorem 2.4 are given in the Appendix. From Theorem 2.4 we observe that in general the bias term asymptotically dominates the other two. It is the unbiasedness of $C_p$ that gives it the asymptotic advantage. In other words, when one compares the asymptotic prediction error for different criteria, essentially the comparison is focused on the bias, and as long as asymptotics is concerned the variability of the criteria does not matter much. Theorem 2.4, therefore, provides an understanding of the gap between finite-sample and asymptotic results regarding selection criteria. Since the asymptotic comparison essentially focuses on the bias and $C_p$ is unbiased, it is not surprising that the high variability of $C_p$ evident in finite-sample studies does not show up in the large-sample considerations. Furthermore, (2.18) and Theorem 2.4 say that for all three selection criteria of interest, GML, $C_p$ and EE, the averaged prediction error $\frac{1}{n}E\|\hat{\mathbf{g}}_{\hat{\lambda}} - \mathbf{g}\|^2$ is of order $O(n^{-4/5})$, an order familiar to many nonparametric problems. Speckman and Sun (2001) studied the asymptotic properties of selection criteria; they showed that GML- and $C_p$- estimated smoothing parameters have the same convergence rate, which, from a different angle, conveys a message similar to Theorem 2.4.

**3. A geometric bridge between the finite-sample and asymptotic results.** In this section, to obtain an intuitive complement to the result of Section 2, we provide a geometric explanation of why the finite-sample variability does not show up in the asymptotics.

3.1. *The geometry of selection criteria.* The fact that $\hat{\lambda}^{(p,q)}$ chooses $\lambda$ as the minimizer of $l_\lambda^{(p,q)}$ implies that $\hat{\lambda}^{(p,q)}$ must satisfy the normal equation $\frac{\partial}{\partial \lambda} l_\lambda^{(p,q)}(\mathbf{z}^{2/q})|_{\lambda=\hat{\lambda}^{(p,q)}} = 0$, which (through simple algebra) can be written as

$$(3.1) \qquad \dot{\boldsymbol{\eta}}_\lambda^{(p,q)\prime}(\mathbf{z}^{2/q} - \boldsymbol{\mu}_\lambda^{(p,q)})|_{\lambda=\hat{\lambda}^{(p,q)}} = 0,$$

where the vector $\dot{\boldsymbol{\eta}}_\lambda^{(p,q)} = (\dot{\eta}_{\lambda 1}^{(p,q)}, \dot{\eta}_{\lambda 2}^{(p,q)}, \ldots, \dot{\eta}_{\lambda n}^{(p,q)})'$, $\dot{\eta}_{\lambda i}^{(p,q)} = -\frac{p}{q\lambda} a_{\lambda i}(c_q b_{\lambda i}^{1/q})^p$, $\boldsymbol{\mu}_\lambda^{(p,q)} = (\mu_{\lambda 1}^{(p,q)}, \mu_{\lambda 2}^{(p,q)}, \ldots, \mu_{\lambda n}^{(p,q)})$ and $\mu_{\lambda i}^{(p,q)} = 1/(c_q b_{\lambda i}^{1/q})$. This normal equation representation suggests a simple geometric interpretation of the $\hat{\lambda}^{(p,q)}$ criterion. For a given observation $\mathbf{z}$, the smoothing parameter is chosen by projecting $z^{2/q}$ onto the line $\{\mu_\lambda^{(p,q)} : \lambda \geq 0\}$ orthogonally to the direction $\dot{\eta}_\lambda^{(p,q)}$. Figure 1 diagrams the geometry two-dimensionally.

In Figure 1 $\mathcal{L}_\lambda^{(p,q)}$ is the hyperplane $\mathcal{L}_\lambda^{(p,q)} = \{\mathbf{z} : (\dot{\boldsymbol{\eta}}_\lambda^{(p,q)})'(\mathbf{z}^{2/q} - \boldsymbol{\mu}_\lambda^{(p,q)}) = 0\}$. Finding the specific hyperplane $\mathcal{L}_\lambda^{(p,q)}$ that passes through $\mathbf{z}^{2/q}$ is equivalent to solving (3.1). It is noteworthy from Figure 1 that different hyperplanes $\mathcal{L}_\lambda^{(p,q)}$ are not parallel, but rather intersect each other, while points on the intersection of two hyperplanes satisfy both normal equations. This phenomenon is termed the *reversal effect* in Efron (2001) and Kou and Efron (2002). Figure 2 provides an illustration, showing one hyperplane $\mathcal{L}^{(p,q)}_{\lambda_0}$ intersecting a nearby hyperplane $\mathcal{L}^{(p,q)}_{\lambda_0+d\lambda}$ (for a small $d\lambda$).

Intuitively, if an observation falls beyond the intersection (i.e., in the reversal region), the selection criterion $\hat{\lambda}^{(p,q)}$ then will have a hard time assigning the smoothing parameter. Furthermore, we observe that for $\hat{\lambda}^{(p,q)}$, if the direction $\dot{\boldsymbol{\eta}}_\lambda^{(p,q)}$ rotates very fast, the reversal region will then be quite large, causing the criterion to have a high chance of encountering observations falling into the reversal region. This reversal effect is the main factor behind $C_p$'s finite-sample unstable behavior, because the $C_p$ orthogonal direction $\dot{\boldsymbol{\eta}}_\lambda^{(2,1)}$ rotates much faster than both the EE direction $\dot{\boldsymbol{\eta}}_\lambda^{(3/2,3/2)}$ and the GML $\dot{\boldsymbol{\eta}}_\lambda^{(1,1)}$ [Kou and Efron (2002)]. It is worth pointing out that the geometry and the reversal effect do not involve asymptotics. Thus finite-sample-wise, the faster rotation of $\dot{\boldsymbol{\eta}}_\lambda^{(2,1)}$ costs $C_p$ much more instability than the EE and GML criteria, undermining its competitiveness.

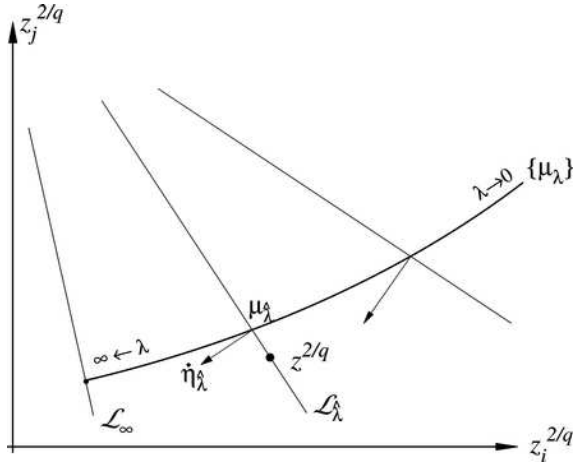

Fig. 1. *The geometry of selection criteria. Two coordinates $z_i^{2/q}$ and $z_j^{2/q}$ ($i < j$) are indicated here.*

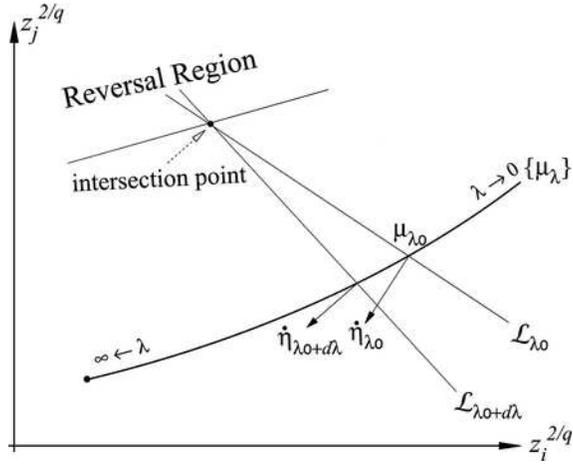

Fig. 2. *Illustration of the reversal effect caused by the rotation of the orthogonal directions.*

3.2. *The evolution of the geometry.* The geometric interpretation naturally suggests we investigate the evolution of the reversal effect (i.e., the geometric instability) as the sample size grows large to bridge the gap between finite- and large-sample results. There are two ways to quantify the geometric instability. First, since the root of instability is the rotation of the orthogonal directions, the curvature of the directions, which captures how fast they rotate, is a measure of the geometric instability. Second, one can investigate the probability that an observation falls into the reversal region, which directly measures how large the reversal effect is.

For the orthogonal direction $\dot{\boldsymbol{\eta}}_\lambda^{(p,q)}$, its statistical curvature [Efron (1975)], which measures the speed of rotation, is defined by

$$\gamma_\lambda = \left(\frac{\det(M_\lambda)}{(\dot{\boldsymbol{\eta}}_\lambda^{(p,q)\prime} V_\lambda \dot{\boldsymbol{\eta}}_\lambda^{(p,q)})^3}\right)^{1/2}$$

$$\text{with } M_\lambda = \begin{pmatrix} \dot{\boldsymbol{\eta}}_\lambda^{(p,q)\prime} V_\lambda \dot{\boldsymbol{\eta}}_\lambda^{(p,q)} & \dot{\boldsymbol{\eta}}_\lambda^{(p,q)\prime} V_\lambda \ddot{\boldsymbol{\eta}}_\lambda^{(p,q)} \\ \dot{\boldsymbol{\eta}}_\lambda^{(p,q)\prime} V_\lambda \ddot{\boldsymbol{\eta}}_\lambda^{(p,q)} & \ddot{\boldsymbol{\eta}}_\lambda^{(p,q)\prime} V_\lambda \ddot{\boldsymbol{\eta}}_\lambda^{(p,q)} \end{pmatrix},$$

where $\ddot{\boldsymbol{\eta}}_\lambda^{(p,q)} = \frac{\partial}{\partial \lambda} \dot{\boldsymbol{\eta}}_\lambda^{(p,q)}$, and the matrix $V_\lambda = \mathrm{diag}(c_q^{-(p+1)} b_{\lambda i}^{-(p+1)/q}/p)$. For the selection criteria class (2.5), Kou and Efron (2002) showed that the squared statistical curvature

$$(3.2) \qquad \gamma_\lambda^2 = \frac{(p+q)^2}{pc_q^{p-1}} \left\{ \frac{\sum_i a_{\lambda i}^4 b_{\lambda i}^{(p-1)/q}}{(\sum_i a_{\lambda i}^2 b_{\lambda i}^{(p-1)/q})^2} - \frac{(\sum_i a_{\lambda i}^3 b_{\lambda i}^{(p-1)/q})^2}{(\sum_i a_{\lambda i}^2 b_{\lambda i}^{(p-1)/q})^3} \right\}.$$

THEOREM 3.1. *The curvature evaluated at the ideal smoothing parameter $\lambda_0$ has the asymptotic order $\gamma_{\lambda_0} = O(n^{-1/10})$.*

PROOF. According to Theorem 2.3, $\gamma_{\lambda_0}^2 = O((\frac{n}{\lambda_0})^{-1/4})$, which is $O(n^{-1/5})$ by (2.15). □

Theorem 3.1 says that, first, for the selection criteria class (2.5), geometrically as the sample size gets larger and larger, the orthogonal directions will rotate more and more slowly, which will make the geometric instability smaller and smaller; second, for different selection criteria, the curvature decreases at the same order.

Next, we consider the probability of an observation falling into the reversal region. Following Kou and Efron (2002), the reversal region (i.e., the region beyond the intersection of different hyperplanes) is defined as

$$\text{reversal region} = \{\mathbf{z} : R_0(\mathbf{z}) < 0\},$$

where the function $R_0(\mathbf{z})$ is given by $R_0(\mathbf{z}) = \ddot{l}^{(p,q)}_{\lambda_0}(\mathbf{z}^{2/q}) - \beta_{\lambda_0} \dot{l}^{(p,q)}_{\lambda_0}(\mathbf{z}^{2/q})$ with $l^{(p,q)}_\lambda$ defined in (2.4), $\dot{l}^{(p,q)}_\lambda = \frac{\partial}{\partial \lambda} l^{(p,q)}_\lambda$, $\ddot{l}^{(p,q)}_\lambda = \frac{\partial^2}{\partial \lambda^2} l^{(p,q)}_\lambda$, and the constant $\beta_{\lambda_0} = -\frac{1}{\lambda_0}[2 - (1 + \frac{p}{q})\frac{\sum_i a^3_{\lambda_0 i} b^{-2/q}_{\lambda_0 i}}{\sum_i a^2_{\lambda_0 i} b^{-2/q}_{\lambda_0 i}}]$.

THEOREM 3.2. *Under mild regularity conditions, the probability that an observation will fall into the reversal region satisfies*

$$P(R_0(\mathbf{z}) < 0) - \Phi(T^{(p,q)}_n) \to 0 \qquad \text{as } n \to \infty,$$

*where $\Phi$ is the standard normal c.d.f. and for all $p \geq q$ the sequence $T^{(p,q)}_n = O(n^{1/10}) < 0$.*

The regularity conditions and proof are deferred to the Appendix. Theorems 3.1 and 3.2 point out that as the sample size $n$ grows large, the reversal effect, which is the source of $C_p$'s instability, decreases at the same rate for all $(p,q)$-estimators and eventually vanishes. This uniform rate is particularly beneficial for $C_p$, because under a finite-sample size, $C_p$ suffers from the reversal effect a lot more than the other criteria, such as GML and EE. Theorems 3.1 and 3.2 thus explain geometrically why the high variability of $C_p$ observed by many authors in finite-sample studies does not hurt it as long as asymptotics is concerned.

**4. A numerical illustration.** In this section through a simulation experiment we will illustrate the connection between finite-sample and asymptotic performances of different selection criteria, focusing on $C_p$, GML and EE. The experiment starts from a small sample size and increases it gradually to exhibit how the performance of different selection criteria evolves as the sample size $n$ grows.

In the simulation the design points $\mathbf{x}$ are $n$ equally spaced points on the $[-1, 1]$ interval, where the sample size $n$ starts at 61, and increases to 121, 241, ..., until 3841. For each value of $n$, 1000 data sets are generated from the curve $f(x) = \sin(\pi(x+1))/(x/2+1)$ shown in Figure 3 with noise level $\sigma = 1$. The $C_p$, GML and EE criteria are applied to the simulated data to choose the smoothing parameter (hence the degrees of freedom), which is subsequently used to estimate the curve.

The bias-variance relationship can be best illustrated by comparing the estimated degrees of freedom (from different selection criteria) with the ideal degrees of freedom $df_0$, since Efron (2001) suggested that the comparison based on degrees of freedom is more sensitive. Figure 4 shows the histograms of $C_p$, GML and EE estimated degrees of freedom under various sample sizes; the vertical bar in each panel represents the ideal degrees of freedom $df_0$.

One can observe from Figure 4 that (i) $C_p$ is roughly unbiased; (ii) as sample size increases, the bias of GML is gradually revealed; (iii) the large spread of $C_p$ estimates points out its high variability even for sample size as large as 3841. The asymptotic results, overlooking the variability, in a certain sense reveal only part of the picture.

Table 1 reports the squared curvature of different selection criteria under various sample sizes; one sees that the curvature of $C_p$ is significantly larger than that of GML or EE, meaning that finite-sample-wise, $C_p$ suffers more from geometric instability. Although the geometric instability (measured by the curvature) becomes smaller and smaller as the sample size gets larger and larger, it decreases quite slowly, indicating that unless the sample size is *very* large, the variability cannot be overlooked (as the asymptotics would do).

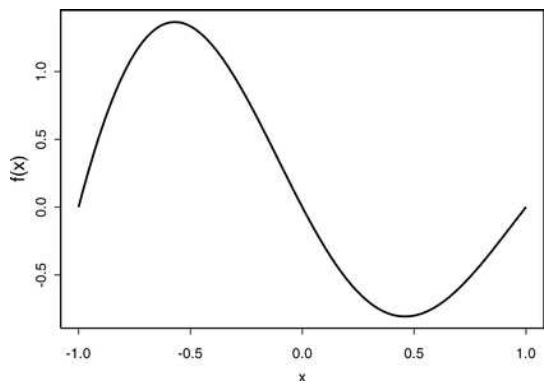

Fig. 3. *The curve used to generate the data.*

Table 1
*The squared curvature of $C_p$, GML and EE*

|       | $n=61$ | $n=121$ | $n=241$ | $n=481$ | $n=961$ | $n=1921$ | $n=3841$ |
|-------|--------|---------|---------|---------|---------|----------|----------|
| $C_p$ | 0.71   | 0.63    | 0.57    | 0.51    | 0.46    | 0.41     | 0.37     |
| GML   | 0.08   | 0.07    | 0.06    | 0.05    | 0.04    | 0.04     | 0.03     |
| EE    | 0.29   | 0.26    | 0.23    | 0.21    | 0.19    | 0.17     | 0.15     |

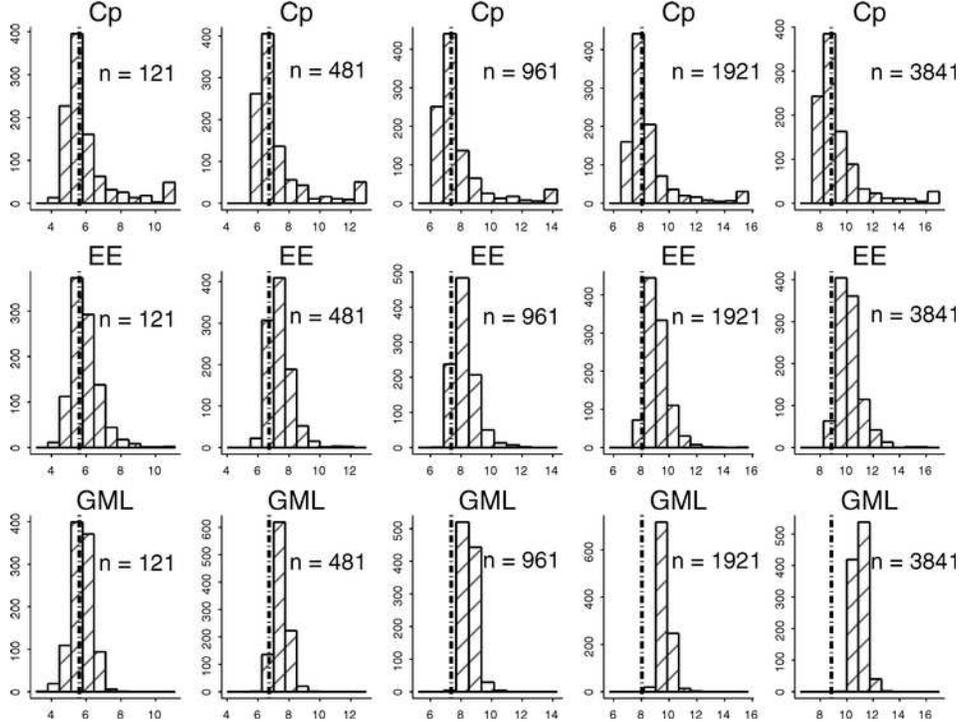

FIG. 4. $C_p$, *GML and EE estimated degrees of freedom. The vertical bar in each panel is the ideal degrees of freedom.*

TABLE 2
*The sample mean and standard deviation of $\|\hat{\mathbf{g}}_{\hat{\lambda}} - \mathbf{g}\|^2$*

|  |  | $n=61$ | $n=121$ | $n=241$ | $n=481$ | $n=961$ | $n=1921$ | $n=3841$ |
|---|---|---|---|---|---|---|---|---|
| $C_p$ | mean | 6.22 | 6.41 | 6.75 | 7.34 | 7.45 | 8.13 | 9.20 |
|  | std dev | 4.81 | 4.54 | 4.42 | 4.42 | 4.25 | 4.33 | 4.91 |
| GML | mean | 5.90 | 5.68 | 5.91 | 6.61 | 7.01 | 7.85 | 9.10 |
|  | std dev | 4.03 | 3.34 | 3.18 | 3.47 | 3.39 | 3.79 | 4.07 |
| EE | mean | 5.89 | 5.78 | 6.10 | 6.73 | 7.03 | 7.78 | 8.86 |
|  | std dev | 4.04 | 3.34 | 3.33 | 3.53 | 3.49 | 3.83 | 4.08 |

Table 2 reports the average value and standard deviation of $\|\hat{\mathbf{g}}_{\hat{\lambda}(p,q)} - \mathbf{g}\|^2$, the squared estimation error, across the data sets. It is interesting to observe that (i) the standard deviation of $C_p$ estimates is larger than that of GML and EE, since geometrically $C_p$ suffers more from the reversal effect than the other two; (ii) for small sample sizes, GML appears to work better than $C_p$ as the asymptotics come in rather slowly; (iii) for reasonable sample sizes from 61 to 3841, as one usually encounters in practice, the EE criterion appears to behave stably well.

Comparing Table 2 with the result of Theorem 2.4, a careful reader might notice that this example itself illustrates the "seeming" gap: For sample size as large as 3841 the asymptotics are still not there. This, again, is due to the fact that although $C_p$'s unbiasedness gives it an asymptotic competitive edge, the asymptotics come in rather slowly, and, therefore, for finite-sample size at hand one cannot neglect the variability, which evidently causes $C_p$ more trouble than the others in Table 2.

**5. Discussion.** This article investigates the connection between finite-sample properties of selection criteria and their asymptotic counterparts, focusing on bridging the gap between the two. Through a bias-variance decomposition of the prediction error, it is shown that in asymptotics bias dominates variability, and thus the large-sample comparison essentially concentrates on bias, and unintentionally overlooks the variability. As the geometry intuitively explains how different selection criteria work, the article also studies the evolution of the geometric instability, the source of $C_p$'s high variability, and shows that although the geometric instability decreases as sample size grows, it decreases very slowly so that for sample sizes one usually encounters in practice, it cannot be neglected. We conclude the article with a few remarks.

REMARK 5.1. *General design points.* We have assumed that the design points $\mathbf{x} = (x_1, \ldots, x_n)$ are equally spaced along a fixed interval. If $\mathbf{x}$ are drawn, instead, from a distribution function $G$ such that $x_i = G^{-1}((2i-1)/n)$, then essentially all the results would remain valid. For example, the conclusion of Theorem 2.3 changes to

$$\sum_{i=3}^{n} a_{\lambda i}^r b_{\lambda i}^s = \frac{1}{4\pi}\left(\int_{\mathcal{X}} g^{1/4}(x)\,dx\right) B\left(r - \frac{1}{4}, s + \frac{1}{4}\right)\left(\frac{n}{\lambda}\right)^{1/4} + o\left(\left(\frac{n}{\lambda}\right)^{1/4}\right),$$

where $g(x)$ is the density of $G$ over the domain $\mathcal{X}$ [Kou (2003)]. Correspondingly, the asymptotic orders that we derived will remain the same (except for longer expressions in the proofs).

REMARK 5.2. *Unknown $\sigma^2$.* To focus on the basic ideas, we implicitly assumed $\sigma^2$ to be known in our analysis. If $\sigma^2$ is unknown, we can replace it with an estimate $\tilde{\sigma}^2$, which changes (2.3) to $\tilde{\mathbf{z}} \equiv \mathbf{U}'y/\tilde{\sigma} = \mathbf{z}(\sigma/\tilde{\sigma})$ and $\tilde{\mathbf{z}}^{2/q} = \mathbf{z}^{2/q}R$, where $R = (\sigma^2/\tilde{\sigma}^2)^{1/q}$, leading to the estimator $\tilde{\lambda}^{(p,q)} = \arg\min_\lambda \{l_\lambda^{(p,q)}(\tilde{\mathbf{z}}^{2/q})\}$, and likewise $\widetilde{df}^{(p,q)}$. If $R \sim (1, \mathrm{var}_R)$ is independent of $\mathbf{z}^{2/q}$, it is easy to see that

(5.1) $\quad \tilde{\mathbf{z}}^{2/q} \sim (E(\mathbf{z}^{2/q}), \mathrm{var}\,\mathbf{z}^{2/q} + \mathrm{var}_R \cdot (E(\mathbf{z}^{2/q})E(\mathbf{z}^{2/q})' + \mathrm{var}\,\mathbf{z}^{2/q})),$

where the notation $X \sim (\alpha, \beta)$ means $X$ has mean $\alpha$ and variance $\beta$. The extra uncertainty of $\sigma^2$ makes the estimate more variable. For example, it can be shown that

$$\frac{\text{var}\{\widetilde{df}^{(p,q)}\}}{\text{var}\{\widehat{df}^{(p,q)}\}} \doteq 1 + \text{var}_R \cdot \left[1 + \frac{(\sum_i a_{\lambda_c i} B_{\lambda_c i}^{p-1}/c_q)^2}{\sum_i a_{\lambda_c i}^2 B_{\lambda_c i}^{2p} \text{var } z_i^{2/q}}\right],$$

which shows the loss of precision in $\widehat{df}^{(p,q)}$ from having to estimate $\sigma^2$. Likewise, our results in Sections 2 and 3 can be modified (at the expense of more complicated calculations) without changing the conclusion. In practice, the estimate $\tilde{\sigma}^2$ can be based on the higher components of $\mathbf{U}'\mathbf{y} \sim (\sigma \mathbf{g}, \sigma^2 \mathbf{I})$, for instance,

$$\tilde{\sigma}^2 = \sum_{i=n-1-M}^{n} (\mathbf{U}'\mathbf{y})_i^2/(M-2),$$

because the assumed smoothness of $\mathbf{f}$ implies that $g_i \doteq 0$ for $i$ large and that $\tilde{\sigma}^2$ and $\mathbf{z}^{2/q}$ are nearly independent, which makes (5.1) valid.

REMARK 5.3. *Higher-order smooth curves.* In Section 2.3, we showed that for general curves EE, $C_p$ and GML gave the same order $O(n^{-4/5})$ for the averaged prediction error $\frac{1}{n}E\|\hat{\mathbf{g}}_{\hat{\lambda}} - \mathbf{g}\|^2$. A reader familiar with the work of Wahba (1985) might sense this as a puzzle, because there it is shown that $C_p$ (GCV) has a faster convergence rate than GML. This seeming conflict actually arises from the difference in the requirements. Wahba (1985) worked on higher-order smooth curves that belong to the null space of the roughness penalty. In our context of cubic smoothing splines they are the curves such that $\int f''(x)^2 dx = 0$, namely, linear lines. In contrast we have assumed $\int f''(x)^2 dx > 0$, and termed them "general curves"; see (2.14).

REMARK 5.4. *Generalizations of $C_p$ and GML.* A number of authors have suggested modifying $C_p$ or GML, including (i) general $C_p$, whose criterion is $C_p(\lambda) = \|\mathbf{y} - \hat{\mathbf{f}}_\lambda\|^2 + 2\omega\sigma^2 \text{tr}(\mathbf{A}_\lambda)$, (ii) general GCV, whose criterion is $GCV(\lambda) = \|\mathbf{y} - \hat{\mathbf{f}}_\lambda\|^2/(1 - \omega \text{ tr}(\mathbf{A}_\lambda)/n)^2$, and (iii) a full Bayesian estimate by putting a prior on the unknown smoothing parameter $\lambda$. Taking $\omega = 1$ in (i) and (ii) results in the classical $C_p$ and GCV. One can also see (through a Taylor expansion) that (i) and (ii) are asymptotically equivalent. Using a number $\omega > 1$ will make the estimate stabler since a heavier roughness penalty is assigned; on the other hand, this will cause the $C_p$ criterion to lose its unbiasedness, since the central smoothing parameter will no longer coincide with the ideal smoothing parameter $\lambda_0$. The finite-sample stability will thus trade off $C_p$'s asymptotic advantage. The full Bayesian approach (iii) is expected to behave even more stably than GML. An interesting open problem is to investigate how large its bias will be and how its geometry, if possible, will evolve as sample size grows.

REMARK 5.5. *Regularity conditions.* All the regularity conditions for the theoretical results, such as Assumptions A.1–A.4 in the Appendix, can be summarized simply as

$$c_q E\{z_i^{2/q}\} \approx 1 + \frac{1}{q} g_i^2,$$

$$\operatorname{var} z_i^{2/q} \approx \operatorname{const} + \operatorname{const} g_i^2,$$

$$E(z_i^{2/q} - E\{z_i^{2/q}\})^3 \approx \operatorname{const} + \operatorname{const} g_i^2.$$

Strict equality holds in the case of $C_p$ and GML, where $q=1$, $c_q=1$:

$$E\{z_i^2\} = 1 + g_i^2,$$

$$\operatorname{var} z_i^2 = 2 + 4g_i^2,$$

$$E(z_i^2 - 1 - g_i^2)^3 = 8 + 24g_i^2,$$

which point out that the conditions are reasonably mild.

## APPENDIX: REGULARITY CONDITIONS AND DETAILED PROOFS

*Regularity conditions for Lemma* 2.2.

ASSUMPTION A.1. $\sum_i a_{\lambda_c i} b_{\lambda_c i}^{p/q} (c_q E\{z_i^{2/q}\} - 1) = O(\sum_i a_{\lambda_c i} b_{\lambda_c i}^{p/q} g_i^2)$.

To see the validity of the assumption, we notice that $q=1$ for $C_p$ and GML, and $\sum_i a_{\lambda_c i} b_{\lambda_c i}^p (c_q E\{z_i^2\} - 1)$ exactly equals $\sum_i a_{\lambda_c i} b_{\lambda_c i}^p g_i^2$. Assumption A.1, hence, clearly holds true for $C_p$ and GML, indicating its mildness. The proof below provides more discussion.

PROOF OF LEMMA 2.2. To prove the lemma, we need the following result of Kou and Efron [(2002), Lemma 1]: For $z_i \sim N(g_i, 1)$, $E(z_i^{2/q}) = \frac{1}{\sqrt{\pi}} 2^{1/q} \Gamma(\frac{1}{q} + \frac{1}{2}) M(-\frac{1}{q}, \frac{1}{2}, -\frac{1}{2} g_i^2)$, where $M(\cdot,\cdot,\cdot)$ is the confluent hypergeometric function (CHF) defined by $M(c,d,z) = 1 + \frac{cz}{d} + \cdots + \frac{(c)_n z^n}{(d)_n n!} + \cdots$, with $(d)_n = d(d+1)\cdots(d+n-1)$. Applying the bounds of CHF [Chapter 13 of Abramowitz and Stegun (1972)]: $1 + \frac{1}{q} g_i^2 - \frac{1}{6q}(1 - \frac{1}{q}) g_i^4 \leq M(-\frac{1}{q}, \frac{1}{2}, -\frac{1}{2} g_i^2) \leq 1 + \frac{1}{q} g_i^2$, one has

(A.1) $$\frac{1}{q} g_i^2 - \frac{1}{6q}\left(1 - \frac{1}{q}\right) g_i^4 \leq c_q E\{z_i^{2/q}\} - 1 \leq \frac{1}{q} g_i^2.$$

The left-hand side of (2.13) is thus bounded above by $\frac{1}{q}\sum_i a_{\lambda_c i} b_{\lambda_c i}^{p/q} g_i^2$, and below by $\frac{1}{q}\sum_i a_{\lambda_c i} b_{\lambda_c i}^{p/q} g_i^2 - \frac{1}{6q}(1-\frac{1}{q})\sum_i a_{\lambda_c i} b_{\lambda_c i}^{p/q} g_i^4$. From (2.14), $\sum_{i=3}^n k_i g_i^2 \leq \frac{1}{\sigma^2}\int_0^1 f''(t)^2\, dt < \infty$, suggesting that for $n$ sufficiently large, the term $\frac{1}{q} g_i^2$ of (A.1) dominates, which again points out that Assumption A.1 is mild. In light of (2.14), $\sum_i a_{\lambda_c i} b_{\lambda_c i}^{p/q} g_i^2 = \lambda_c \sum_i a_{\lambda_c i}^2 b_{\lambda_c i}^{p/q-1}(k_i g_i^2) = O(\lambda_c)$, for $p \geq q$, which according to Assumption A.1 implies that

$$\sum_i a_{\lambda_c i} b_{\lambda_c i}^{p/q} (c_q E\{z_i^{2/q}\} - 1) = O(\lambda_c) \qquad \text{for } p \geq q.$$

$\square$

To prove Theorem 2.4, we need the following approximation.

LEMMA A.1.

(A.2)
$$\begin{aligned}
E\|\hat{\mathbf{g}}_{\hat{\lambda}} - \hat{\mathbf{g}}_{\lambda_c}\|^2 \\
&\doteq \frac{c_q^2}{Q_{\lambda_c}^2(E\{\mathbf{z}^{2/q}\})} \\
&\quad \times \left\{\left(\sum_i a_{\lambda_c i}^2 b_{\lambda_c i}^2 (g_i^2 + 1)\right)\left(\sum_i a_{\lambda_c i}^2 b_{\lambda_c i}^{2p/q} \operatorname{var} z_i^{2/q}\right) \right. \\
&\qquad \left. + \sum_i a_{\lambda_c i}^4 b_{\lambda_c i}^{2+2p/q} E[(z_i^2 - g_i^2 - 1)(z_i^{2/q} - E\{z_i^{2/q}\})^2]\right\},
\end{aligned}$$

(A.3)
$$\begin{aligned}
E(\hat{\mathbf{g}}_{\lambda_c} - \mathbf{g})'(\hat{\mathbf{g}}_{\hat{\lambda}} - \hat{\mathbf{g}}_{\lambda_c}) \\
&\doteq \frac{c_q}{Q_{\lambda_c}(E\{z^{2/q}\})} \\
&\quad \times \sum_i a_{\lambda_c i}^2 b_{\lambda_c i}^{1+p/q}(a_{\lambda_c i}\operatorname{cov}(z_i^2, z_i^{2/q}) - g_i \operatorname{cov}(z_i, z_i^{2/q})),
\end{aligned}$$

where the function $Q_\lambda(\mathbf{u})$ is defined by

(A.4) $\quad Q_\lambda(\mathbf{u}) = \sum_i a_{\lambda i} b_{\lambda i}^{(p-1)/q}\left\{\frac{1}{q}a_{\lambda i} + \left[\left(1 + \frac{p}{q}\right)a_{\lambda i} - 2\right](c_q b_{\lambda i}^{1/q} u_i - 1)\right\}.$

DERIVATION OF LEMMA A.1. Since $\hat{\lambda}$ by definition is a function of $\mathbf{u} = \mathbf{z}^{2/q}$, and $\lambda_c$ is a function of $E\{\mathbf{z}^{2/q}\}$, applying a Taylor expansion on $a_{\hat{\lambda} i} - a_{\lambda_c i}$, we obtain

$$a_{\hat{\lambda} i} - a_{\lambda_c i} \doteq -\frac{a_{\lambda_c i} b_{\lambda_c i}}{\lambda_c} \cdot \sum_j \left.\frac{\partial \hat{\lambda}}{\partial u_j}\right|_{\mathbf{u}=E\{\mathbf{z}^{2/q}\}} (z_j^{2/q} - E\{z_j^{2/q}\}).$$

Some algebra, after applying the implicit function calculation to the definition (2.5) of $\hat{\lambda}$ or equivalently to the normal equation (3.1), yields $\frac{\partial \hat{\lambda}}{\partial u_j}|_{\mathbf{u}=E\{\mathbf{z}^{2/q}\}} = \frac{-\lambda_c c_q a_{\lambda j} b_{\lambda j}^{p/q}}{Q_{\lambda_c}(E\{z^{2/q}\})}$, which then gives

$$(A.5) \quad \hat{g}_{\hat{\lambda}i} - \hat{g}_{\lambda_c i} = (a_{\hat{\lambda}i} - a_{\lambda_c i})z_i \doteq \frac{c_q a_{\lambda_c i} b_{\lambda_c i} z_i}{Q_{\lambda_c}(E\{z^{2/q}\})} \sum_j a_{\lambda_c j} b_{\lambda_c j}^{p/q}(z_j^{2/q} - E\{z_j^{2/q}\}).$$

The fact that the $z_i$'s are independent of each other implies

$$E(\hat{g}_{\hat{\lambda}i} - \hat{g}_{\lambda_c i})^2$$
$$\doteq \frac{c_q^2 a_{\lambda_c i}^2 b_{\lambda_c i}^2}{Q_{\lambda_c}^2(E\{z^{2/q}\})} \bigg\{ (g_i^2 + 1) \sum_j a_{\lambda_c j}^2 b_{\lambda_c j}^{2p/q} \operatorname{var} z_j^{2/q}$$
$$+ a_{\lambda_c i}^2 b_{\lambda_c i}^{2p/q} E[(z_i^2 - g_i^2 - 1)(z_i^{2/q} - E\{z_i^{2/q}\})^2] \bigg\}.$$

Summing over $i$ yields the approximation

$$E\|\hat{\mathbf{g}}_{\hat{\lambda}} - \hat{\mathbf{g}}_{\lambda_c}\|^2$$
$$\doteq \frac{c_q^2}{Q_{\lambda_c}^2(E\{\mathbf{z}^{2/q}\})} \bigg\{ \bigg(\sum_i a_{\lambda_c i}^2 b_{\lambda_c i}^2 (g_i^2 + 1)\bigg) \bigg(\sum_i a_{\lambda_c i}^2 b_{\lambda_c i}^{2p/q} \operatorname{var} z_i^{2/q}\bigg)$$
$$+ \sum_i a_{\lambda_c i}^4 b_{\lambda_c i}^{2+2p/q} E[(z_i^2 - g_i^2 - 1)(z_i^{2/q} - E\{z_i^{2/q}\})^2] \bigg\}.$$

The approximation of $E(\hat{\mathbf{g}}_{\lambda_c} - \mathbf{g})'(\hat{\mathbf{g}}_{\hat{\lambda}} - \hat{\mathbf{g}}_{\lambda_c})$ can be obtained in a similar way. □

Before proving Theorem 2.4, we state its regularity conditions. Theorem 2.4 needs the following assumptions, in addition to Assumption A.1.

ASSUMPTION A.2.
$$\sum_i a_{\lambda_c i} b_{\lambda_c i}^{p/q} \bigg[\bigg(1 + \frac{p}{q}\bigg) a_{\lambda_c i} - 2\bigg] (c_q E\{z_i^{2/q}\} - 1)$$
$$= O\bigg(\sum_i a_{\lambda_c i} b_{\lambda_c i}^{p/q} \bigg[\bigg(1 + \frac{p}{q}\bigg) a_{\lambda_c i} - 2\bigg] g_i^2 \bigg),$$
$$\sum_i a_{\lambda_c i}^2 b_{\lambda_c i}^{2p/q} \operatorname{var} z_i^{2/q}$$
$$= O\bigg(\max\bigg(\sum_i a_{\lambda_c i}^2 b_{\lambda_c i}^{2p/q}, \sum_i a_{\lambda_c i}^2 b_{\lambda_c i}^{2p/q} g_i^2\bigg)\bigg).$$

ASSUMPTION A.3.   $\sum_i (a_{\lambda_c i}^2 b_{\lambda_c i}^{1+p/q})^l E\{z_i^m z_i^{2n/q}\} = O(\max(\sum_i (a_{\lambda_c i}^2 b_{\lambda_c i}^{1+p/q})^l, \sum_i (a_{\lambda_c i}^2 b_{\lambda_c i}^{1+p/q})^l g_i^2))$, for $l, m, n \in \{1, 2\}$.

Like Assumption A.2, these two assumptions are exactly true for GML and $C_p$, since $E\{z_i^2\} = 1 + g_i^2$, and $\text{var}(z_i^2) = 2 + 4g_i^2$. In general, a Taylor expansion on the CHF can show for $q \geq 1$,

$$\text{(A.6)} \qquad \text{var } z_i^{2/q} = \text{const } 1 + \text{const } 2 \cdot g_i^2 + O(g_i^4),$$

which suggests that the assumptions are mild.

PROOF OF THEOREM 2.4.   Write

$$\text{Term A} = \left(\sum_i a_{\lambda_c i}^2 b_{\lambda_c i}^2 (g_i^2 + 1)\right)\left(\sum_i a_{\lambda_c i}^2 b_{\lambda_c i}^{2p/q} \text{var } z_i^{2/q}\right),$$

$$\text{Term B} = \sum_i a_{\lambda_c i}^4 b_{\lambda_c i}^{2+2p/q} E[(z_i^2 - g_i^2 - 1)(z_i^{2/q} - E\{z_i^{2/q}\})^2];$$

then approximation (A.2) becomes

$$\text{(A.7)} \qquad E\|\hat{\mathbf{g}}_{\hat{\lambda}} - \hat{\mathbf{g}}_{\lambda_c}\|^2 \doteq \frac{c_q^2}{Q_{\lambda_c}^2(E\{\mathbf{z}^{2/q}\})}(\text{Term A} + \text{Term B}).$$

For Term A, note that according to Assumption A.2 the order of $\sum_i a_{\lambda_c i}^2 b_{\lambda_c i}^{2p/q} \times \text{var } z_i^{2/q}$ is the maximum of $\sum_i a_{\lambda_c i}^2 b_{\lambda_c i}^{2p/q}$ and $\sum_i a_{\lambda_c i}^2 b_{\lambda_c i}^{2p/q} g_i^2$. But $\sum_i a_{\lambda_c i}^2 b_{\lambda_c i}^{2p/q} = O((\frac{n}{\lambda_c})^{1/4}) = O(n^{1/5})$ by Theorem 2.3, and $\sum_i a_{\lambda_c i}^2 b_{\lambda_c i}^{2p/q} g_i^2 = O(\lambda_c) = O(n^{1/5})$. So $\sum_i a_{\lambda_c i}^2 b_{\lambda_c i}^{2p/q} \text{var } z_i^{2/q} = O(n^{1/5})$. Next observe that $\sum_i a_{\lambda_c i}^2 b_{\lambda_c i}^2 (g_i^2 + 1) = \sum_i a_{\lambda_c i}^2 b_{\lambda_c i}^2 g_i^2 + \sum_i a_{\lambda_c i}^2 b_{\lambda_c i}^2$; the first term is equal to $\lambda_c (\sum_i a_{\lambda_c i}^3 b_{\lambda_c i}(k_i g_i^2)) = O(\lambda_c) = O(n^{1/5})$; the second term is of order $O((\frac{n}{\lambda_c})^{1/4}) = O(n^{1/5})$. Therefore, Term A $= O(n^{1/5} \cdot n^{1/5}) = O(n^{2/5})$.

For Term B, a Taylor expansion on the CHF gives

$$\text{(A.8)} \quad E[(z_i^2 - g_i^2 - 1)(z_i^{2/q} - E\{z_i^{2/q}\})^2] = \text{const} + \text{const } g_i^2 + O(g_i^4),$$

which, together with Assumption A.3, implies that the order of Term B is the maximum of $O(\lambda_c) = O(n^{1/5})$ and $O((\frac{n}{\lambda_c})^{1/4}) = O(n^{1/5})$. Thus Term B $= O(n^{1/5})$.

Using Assumption A.2, the denominator in (A.7)

$$Q_{\lambda_c}(E\{\mathbf{z}^{2/q}\}) = \sum_i a_{\lambda_c i} b_{\lambda_c i}^{(p-1)/q}\left\{\frac{1}{q}a_{\lambda_c i} + (b_{\lambda_c i}^{1/q} - 1)\left[\left(1 + \frac{p}{q}\right)a_{\lambda_c i} - 2\right]\right\}$$

$$+ \sum_i a_{\lambda_c i} b_{\lambda_c i}^{p/q}\left[\left(1 + \frac{p}{q}\right)a_{\lambda_c i} - 2\right](c_q E\{z_i^{2/q}\} - 1)$$

(A.9)
$$= O\left(\left(\frac{n}{\lambda_c}\right)^{1/5}\right) + O(\lambda_c)$$
$$= O(n^{1/5}).$$

Plugging (A.9) and the orders of Term A and Term B into (A.7) yields $E\|\hat{\mathbf{g}}_{\hat{\lambda}} - \hat{\mathbf{g}}_{\lambda_c}\|^2 = O(1)$.

For the covariance term $E(\hat{\mathbf{g}}_{\lambda_c} - \mathbf{g})'(\hat{\mathbf{g}}_{\hat{\lambda}} - \hat{\mathbf{g}}_{\lambda_c})$, we can write

$$E(\hat{\mathbf{g}}_{\lambda_c} - \mathbf{g})'(\hat{\mathbf{g}}_{\hat{\lambda}} - \hat{\mathbf{g}}_{\lambda_c}) \doteq \frac{c_q}{Q_{\lambda_c}(E\{z^{2/q}\})}(\text{Term C} + \text{Term D}),$$

$$\text{Term C} = \sum_i a_{\lambda_c i}^2 b_{\lambda_c i}^{1+p/q} \operatorname{cov}(z_i^2, z_i^{2/q}),$$

$$\text{Term D} = -\sum_i a_{\lambda_c i}^2 b_{\lambda_c i}^{1+p/q} g_i \operatorname{cov}(z_i, z_i^{2/q}).$$

Applying Assumption A.3 and the facts that

$$\operatorname{cov}(z_i^2, z_i^{2/q}) = \text{const} + \text{const } g_i^2 + O(g_i^4),$$
$$g_i \operatorname{cov}(z_i, z_i^{2/q}) = \text{const } g_i^2 + O(g_i^4),$$

which can be derived similarly to (A.8), it can be shown that

$$\text{Term C} = O(n^{1/5}), \qquad \text{Term D} = O(n^{1/5}),$$

which finally gives $E(\hat{\mathbf{g}}_{\lambda_c} - \mathbf{g})'(\hat{\mathbf{g}}_{\hat{\lambda}} - \hat{\mathbf{g}}_{\lambda_c}) = O(1)$. □

*Regularity conditions for Theorem* 3.2.

ASSUMPTION A.4.

$$\sum_i [a_{\lambda_0 i}^l b_{\lambda_0 i}^{p/q}(c_q E\{z_i^{2/q}\} - 1)] = O\left(\sum_i a_{\lambda_0 i}^l b_{\lambda_0 i}^{p/q} g_i^2\right) \qquad \text{for } l = 1, 2,$$

$$\sum_i [a_{\lambda_0 i}^l b_{\lambda_0 i}^{2p/q} \operatorname{var} z_i^{2/q}]$$
$$= O\left(\max\left(\sum_i a_{\lambda_0 i}^l b_{\lambda_0 i}^{2p/q}, \sum_i a_{\lambda_0 i}^l b_{\lambda_0 i}^{2p/q} g_i^2\right)\right) \qquad \text{for } l = 2, 3, 4.$$

Like the previous three assumptions, Assumption A.4 is exact for GML and $C_p$. For general criteria in the class (2.5), the facts (A.6), (A.1) and $E(z_i^{2/q} - E\{z_i^{2/q}\})^3 = \text{const} + \text{const } g_i^2 + O(g_i^4)$ suggest that Assumption A.4 is reasonably mild.

PROOF OF THEOREM 3.2. Let $M(R_0)$ and $V(R_0)$ denote the mean and variance of $R_0(\mathbf{z})$. Kou and Efron (2002) showed that

$$M(R_0) = \frac{p}{q^2}(p+q)c_q^{p-1}$$
$$\times \left\{ \frac{1}{p+q}\left(\sum_i a_{\lambda_0 i}^2 b_{\lambda_0 i}^{(p-1)/q}\right) \right.$$
$$\left. + \sum_i \left[ a_{\lambda_0 i} b_{\lambda_0 i}^{(p-1)/q}\left(a_{\lambda_0 i} - \frac{\sum_i a_{\lambda_0 i}^3 b_{\lambda_0 i}^{-2/q}}{\sum_i a_{\lambda_0 i}^2 b_{\lambda_0 i}^{-2/q}}\right)(c_q b_{\lambda_0 i}^{1/q} E\{z_i^{2/q}\} - 1) \right] \right\},$$

$$V(R_0) = \frac{p^2}{q^4}(p+q)^2 c_q^{2p} \sum_i \left[ a_{\lambda_0 i}^2 b_{\lambda_0 i}^{2p/q}\left(a_{\lambda_0 i} - \frac{\sum_i a_{\lambda_0 i}^3 b_{\lambda_0 i}^{-2/q}}{\sum_i a_{\lambda_0 i}^2 b_{\lambda_0 i}^{-2/q}}\right)^2 \operatorname{var} z_i^{2/q} \right].$$

Using the Berry–Esseen theorem [Feller (1971), page 521], we have

$$P(R_0(\mathbf{z}) < 0) - \Phi(M(R_0)/\sqrt{V(R_0)}) \to 0 \quad \text{as } n \to \infty.$$

Note that we can write $\frac{M(R_0)}{\sqrt{V(R_0)}} = \frac{\text{Term 1} + \text{Term 2}}{c_q(\text{Term 3})^{1/2}}$, where

$$\text{Term 1} = \frac{1}{p+q}\left(\sum_i a_{\lambda_0 i}^2 b_{\lambda_0 i}^{(p-1)/q}\right)$$
$$+ \sum_i \left[ a_{\lambda_0 i} b_{\lambda_0 i}^{(p-1)/q}\left(a_{\lambda_0 i} - \frac{\sum_i a_{\lambda_0 i}^3 b_{\lambda_0 i}^{-2/q}}{\sum_i a_{\lambda_0 i}^2 b_{\lambda_0 i}^{-2/q}}\right)(b_{\lambda_0 i}^{1/q} - 1) \right],$$

$$\text{Term 2} = \sum_i \left[ a_{\lambda_0 i} b_{\lambda_0 i}^{p/q}\left(a_{\lambda_0 i} - \frac{\sum_i a_{\lambda_0 i}^3 b_{\lambda_0 i}^{-2/q}}{\sum_i a_{\lambda_0 i}^2 b_{\lambda_0 i}^{-2/q}}\right)(c_q E\{z_i^{2/q}\} - 1) \right]$$

and

$$\text{Term 3} = \sum_i \left[ a_{\lambda_0 i}^2 b_{\lambda_0 i}^{2p/q}\left(a_{\lambda_0 i} - \frac{\sum_i a_{\lambda_0 i}^3 b_{\lambda_0 i}^{-2/q}}{\sum_i a_{\lambda_0 i}^2 b_{\lambda_0 i}^{-2/q}}\right)^2 \operatorname{var} z_i^{2/q} \right].$$

To obtain the order of Term 1, we need another result from Kou (2003): Suppose $\frac{n}{\lambda} \to \infty$; then for all $r > \frac{1}{4}$ and $s < -\frac{1}{4}$, $\sum_{i=3}^n a_{\lambda i}^r b_{\lambda i}^s = O((\frac{n}{\lambda})^{-s})$. This result and Theorem 2.3 imply

(A.10) $$\text{Term 1} = O\left(\left(\frac{n}{\lambda_0}\right)^{1/4}\right) = O(n^{1/5}).$$

To obtain the order of Term 2, we note that by Assumption A.4 and (2.14),

$$\text{Term 2} = O\left(\lambda_0 \sum_i \left[ a_{\lambda_0 i}^2 b_{\lambda_0 i}^{p/q-1}\left(a_{\lambda_0 i} - \frac{\sum_i a_{\lambda_0 i}^3 b_{\lambda_0 i}^{-2/q}}{\sum_i a_{\lambda_0 i}^2 b_{\lambda_0 i}^{-2/q}}\right)(k_i g_i^2) \right]\right)$$

(A.11)
$$= O(\lambda_0) = O(n^{1/5}) \qquad \text{for all } p \geq q.$$

For Term 3, since

$$\sum_i \left[ a_{\lambda_0 i}^2 b_{\lambda_0 i}^{2p/q} \left( a_{\lambda_0 i} - \frac{\sum_i a_{\lambda_0 i}^3 b_{\lambda_0 i}^{-2/q}}{\sum_i a_{\lambda_0 i}^2 b_{\lambda_0 i}^{-2/q}} \right)^2 \right] = O\left( \left( \frac{n}{\lambda_0} \right)^{1/4} \right) = O(n^{1/5})$$

and

$$\sum_i \left[ a_{\lambda_0 i}^2 b_{\lambda_0 i}^{2p/q} \left( a_{\lambda_0 i} - \frac{\sum_i a_{\lambda_0 i}^3 b_{\lambda_0 i}^{-2/q}}{\sum_i a_{\lambda_0 i}^2 b_{\lambda_0 i}^{-2/q}} \right)^2 g_i^2 \right]$$

$$= \lambda_0 \sum_i \left[ a_{\lambda_0 i}^3 b_{\lambda_0 i}^{2p/q-1} \left( a_{\lambda_0 i} - \frac{\sum_i a_{\lambda_0 i}^3 b_{\lambda_0 i}^{-2/q}}{\sum_i a_{\lambda_0 i}^2 b_{\lambda_0 i}^{-2/q}} \right)^2 (k_i g_i^2) \right]$$

$$= O(\lambda_0) = O(n^{1/5}) \qquad \text{for all } p \geq q,$$

using Assumption A.4 we have

(A.12) $\qquad \text{Term } 3 = O(n^{1/5}) \qquad \text{for all } p \geq q.$

Combining (A.10)–(A.12) finally yields

$$T_n^{(p,q)} = \frac{M(R_0)}{\sqrt{V(R_0)}} = O\left( \frac{n^{1/5}}{n^{1/10}} \right) = O(n^{1/10}) < 0 \qquad \text{for all } p \geq q. \qquad \square$$

**Acknowledgments.** The author is grateful to Professor Bradley Efron for helpful discussions. The author also thanks the Editor, the Associate Editor and two referees for constructive suggestions, which much improved the presentation of the paper.

Department of Statistics  
Science Center 6th Floor  
Harvard University  
Cambridge, Massachusetts 02138  
USA  
e-mail: kou@stat.harvard.edu